\documentclass[11pt]{proc-l}
\usepackage{amsmath,amssymb,amsbsy,amsfonts,amsthm,latexsym,
            amsopn,amstext,amsxtra,euscript,amscd,amsthm}

\newtheorem{thm}{Theorem}
\newtheorem{theorem}[thm]{Theorem}

\def\\{\cr}
\def\({\left(}
\def\){\right)}
\def\[{\left[}
\def\]{\right]}
\def\<{\langle}
\def\>{\rangle}

\begin{document}

\title{Cullen numbers with the Lehmer property}

\author{\sc Jos\'e~Mar\'{\i}a~Grau~Ribas}
\address{Departamento de M\'atematicas,\\
Universidad de Oviedo\\
 Avda. Calvo Sotelo, s/n, 33007 Oviedo, Spain}
\email{grau@uniovi.es}

\author{Florian~Luca} 
\address{Instituto de Matem{\'a}ticas,\\
Universidad Nacional Autonoma de M{\'e}xico,\\
C.P. 58089, Morelia, Michoac{\'a}n, M{\'e}xico}
\email{fluca@matmor.unam.mx}

\maketitle

\begin{abstract} 
Here, we show that there is no positive integer $n$ such that the $n$th Cullen number $C_n=n2^n+1$ has the property that it is composite but $\phi(C_n)\mid C_n-1$.
\end{abstract}

\section{Introduction}

A {\it Cullen number} is a number of the form $C_n=n2^n+1$ for some $n\ge 1$. They attracted attention of researchers since it seems that it is hard to find primes of this form. Indeed, Hooley \cite{Ho} showed that for most $n$ the number $C_n$ is composite. For more about testing $C_n$ for primality, see \cite{BF} and \cite{GO}. For an integer $a>1$, a {\it pseudoprime} to base $a$ is a compositive positive integer $m$ such that $a^m\equiv a\pmod m$. Pseudoprime Cullen numbers have also been studied. For example, in 
\cite{LuSh} it is shown that for most $n$, $C_n$ is not a base $a$-pseudoprime. Some computer searchers up to several millions did not turn up any pseudo-prime $C_n$ to any base. Thus, it would seem that Cullen numbers which are pseudoprimes are very scarce.  A {\it Carmichael number} is a positive integer $m$ which is a base $a$ pseudoprime for any $a$. A composite integer $m$ is called a {\it Lehmer number} if $\phi(m)\mid m-1$, where $\phi(m)$ is the Euler function of $m$. Lehmer numbers are Carmichael numbers; hence, pseudoprimes in every base. No Lehmer number is known, although it is known that there are no Lehmer numbers in certain sequences, such as the Fibonacci sequence (see \cite{Lu1}), or the sequence of repunits in base $g$ for any $g\in [2,1000]$ (see \cite{CiLu}). For other results on Lehmer numbers, see \cite{BaGu}, \cite{BaLu}, \cite{LuPo}, \cite{Po1}, \cite{Po2}.

\medskip

Our result here is that there is no Cullen number with the Lehmer property.  Hence, if $\phi(C_n)\mid C_n-1$, then $C_n$ is prime.

\begin{theorem}
\label{thm:main}
Let $C_n$ be the $n$th Cullen number. If $\phi(C_n) \mid C_n-1$, then $C_n$ is prime.
\end{theorem}

\section{Proof of Theorem \ref{thm:main}}

Assume that $n\ge 30$, that $\phi(C_n)\mid C_n-1$, but that $C_n$ is not prime. Then $C_n$ is square-free. Write
$$
C_n=\prod_{i=1}^k p_i.
$$
So,
$$
\prod_{i=1}^k (p_i-1)\mid n2^n.
$$
Write $n=2^{\alpha}n_1$, where $n_1$ is odd. Then 
$C_n=n_12^{n_2}+1$, where $n_2:=\alpha+n$. Let $p$ be any  prime factor of $C_n$. Since $p-1\mid C_n-1$, it follows that $p=m_p2^{n_p}+1$ for some odd divisor $m_p$ of $n$ and some $n_p$ with 
$$
n_p\le n_2=n+\alpha\le n+\frac{\log n}{\log 2}.
$$ 
Let us first show that in fact $n_p\le n$. Assume that $n_p>n$. Then, 
\begin{equation}
\label{eq:11}
C_n=n2^n+1=p\lambda,
\end{equation}
for some positive integer $\lambda$, where $p\ge 2^{n+1}+1$. Observe that $\lambda>1$ because $C_n$ is not prime. Now 
$$
\lambda=\frac{C_n}{p}\le \frac{n2^n+1}{2^{n+1}+1}<n.
$$
Reducing equation \eqref{eq:11} modulo $2^n$, we get that $2^n\mid \lambda-1$, so $2^{n}\le \lambda-1<n$, which is false for any $n>1$. Hence, $n_p\le n$.

\medskip

Next we look at $m_p$. If $m_p=1$, then $p=2^{n_p}+1$ is a Fermat prime. Hence, $n_p=2^{\gamma_p}$ for some nonnegative integer $\gamma$. Since $2^{\gamma_p}=n_p\le n$, we get that $\gamma_p< (\log n)/(\log 2)$. Hence, the prime $p$ can take at most 
$1+(\log  n)/(\log 2)$ values. Next, observe that since
\begin{equation}
\label{eq:mp}
\prod_{p\mid C_n} m_p\mid n,
\end{equation}
it follows that the number of prime factors $p$ of $C_n$ such that $m_p>1$ is $\le (\log n)/(\log 3)$. Hence, we arrived at the bound
\begin{equation}
\label{eq:k}
k<1+\frac{\log n}{\log 2}+\frac{\log n}{\log 3}<1+2.4\log n.
\end{equation}

\medskip

We next bound $n_p$. Put $N:=\lfloor {\sqrt{n/\log n}}\rfloor$, and consider pairs $(a,b)$ of integers in 
$\{0,1,\ldots,N\}$. There are $(N+1)^2>n/\log n$ such pairs. For each such pair, consider the expression $L(a,b):=an+bn_p\in [0,2n^{3/2}/(\log n)^{1/2}]$. Thus, there exist two pairs $(a,b)\ne (a_1,b_1)$ such that 
$$
|(a-a_1)n+(b-b_1)n_p|=|L(a,b)-L(a_1,b_1)|\le \frac{2n^{3/2}/(\log n)^{1/2}}{n/\log n-1}<3(n\log n)^{1/2}.
$$
Put $u:=a-a_1,~v:=b-b_1$. Then $(u,v)\ne (0,0)$ and 
$$
|un+vn_p|<3(n\log n)^{1/2}.
$$ 
We may also assume that $u$ and $v$ are coprime, for if not, we replace the pair $(u,v)$ by the pair $(u_1,v_1)$, where $d:=\gcd(u,v),~u_1:=u/d,~v_1:=v/d$, and the properties that $\max\{|u_1|,|v_1|\}\le (n/\log n)^{1/2}$ and $|u_1n+v_1 n_p|<3(n\log n)^{1/2}$ are still fulfilled. Finally, up to replacing the pair $(u,v)$ by the pair $(-u,-v)$, we may assume that $u\ge 0$. 

Now consider the congruences $n2^n\equiv -1\pmod p$ and $m_p2^{n_p}\equiv -1\pmod p$. Observe that $2,~n,~m_p$ are all three coprime to $p$. Raise the first congruence to $u$ and the second to $v$ and multiply them to get
$$
n^u m_p^v2^{nu+n_pv}\equiv (-1)^{u+v}\pmod p.
$$
Hence, $p$ divides the numerator of the rational number 
\begin{equation}
\label{eq:xy}
A:=n^u m_p^v 2^{nu+n_p v}-(-1)^{u+v}.
\end{equation}
Let us show that $A\ne 0$. Assume that $A=0$. Recall that $C_n=n_12^{n_2}+1$. Thus, expression \eqref{eq:xy} is 
$$
A=n_1^u m_p^v 2^{(n+\alpha)u+n_pv}-(-1)^{u+v}=0.
$$
Then 
$n_1^u m_p^v=1$, $(n+\alpha)u+vn_p=0$, and $u+v$ is even. Since $u\ge 0$, it follows that $v\le 0$. Put $w:=-v$, so $w\ge 0$. There exists a positive integer  $\rho$ which is odd such that $n_1=\rho^{w}$ and $m_p=\rho^u$. Since $u$ and $v$ are coprime and $u+v$ is even, it follows that $u$ and $v$ are both odd. Hence, $w$ is also odd. Also, since $m_p$ divides $n_1$, it follows that $u\le w$. We now get
$$
(2^{\alpha} \rho^{w}+\alpha) u-w n_p=0,
$$
so
$$
\frac{u}{n_p}=\frac{w}{2^{\alpha} \rho^w+\alpha}.
$$
The left--hand side is $\ge u/n=u/(2^{\alpha}\rho^u)$, because $n_p\le n=2^{\alpha}\rho^u$. Hence, we get that
$$
\frac{u}{2^{\alpha}\rho^u}\le \frac{u}{n_p}=\frac{w}{2^{\alpha} \rho^w+\alpha}\qquad {\text{\rm leading~to}}\qquad \frac{u}{\rho^{u}}\le 
\frac{w}{\rho^{w}+(\alpha/2^{\alpha})}\le \frac{w}{\rho^{w}}.
$$
For $\rho\ge 3$, the function $s\mapsto s/\rho^s$ is decreasing for $s\ge 0$, so the above inequality together with the fact that $u\le w$ implies that  $u=w$ (so both are $1$ because they are coprime), and that all the intermediary inequalities are also equalities. This means that $u=w=1$, 
$\alpha=0$ and $n=n_p$, but all this is possible only when $C_n=p$, which is not allowed. If $\rho=1$, we then get that $n_1=1$, so 
every prime factor $p$ of $C_n$ is a Fermat prime. Hence, we get
$$
C_n=2^{n_2}+1=\prod_{i=1}^{k} (2^{2^{\gamma_{p_i}}}+1)=\sum_{I\subseteq \{1,\ldots,k\}} 2^{\sum_{i\in  I} 2^{\gamma_{p_i}}},
$$
and $k\ge 2$, but this is impossible by the unicity of the binary expansion of $C_n$.

Thus, it is not possible for the expression $A$ shown at \eqref{eq:xy} to be zero. 

The size of the numerator of $A$ is at most
\begin{eqnarray*}
& & 2^{1+|nu+n_pv|} n^{u} m_p^{|v|}\le 2^{1+3(n\log n)^{1/2}} n^{2(n/\log n)^{1/2}}\\
& < & 2^{1+3(n\log n)^{1/2}+(2/\log 2)(n\log n)^{1/2}}<2^{6(n\log n)^{1/2}}.
\end{eqnarray*}
In the above chain of inequalities, we used the fact that $3+2/\log 2<5.9$, together with the fact that $(n\log n)^{1/2}>10$ for $n\ge 30$.
Thus, for $n\ge 30$, we have that the inequality
\begin{equation}
\label{eq:p}
p<2^{6(n\log n)^{1/2}}
\end{equation}
holds for all prime factors $p$ of $C_n$.

Thus, we get the inequality
\begin{equation*}
2^n<C_n=\prod_{i=1}^k p_i<\prod_{i=1}^k 2^{6(n\log n)^{1/2}}=2^{6k(n\log n)^{1/2}},
\end{equation*}
leading to
\begin{equation}
\label{eq:low}
k>\frac{n^{1/2}}{6(\log n)^{1/2}}.
\end{equation}

Comparing estimates \eqref{eq:k} and \eqref{eq:low}, we get
$$
\frac{n^{1/2}}{6(\log n)^{1/2}}<1+2.4\log n,
$$
implying $n<6\times 10^5$. 

It remains to lower this bound. We first lower it to $n<93000$. Indeed, first note that since $n<6\times 10^5$, it follows that if $p=F_{\gamma}=2^{2^{\gamma}}+1$ is a Fermat prime dividing $C_n$, then $\gamma\le 18$. The only such Fermat primes are for $\gamma\in \{0,1,2,3,4\}$. Furthermore, $(\log n)/(\log 3)\le \log(6\times 10^5)/(\log 3)=12.1104\ldots$ Hence, $k\le 5+12=17$. It then follows, by equation \eqref{eq:low}, that
$$
\frac{n^{1/2}}{6(\log n)^{1/2}}<17,
$$
so $n<122000$. But then $(\log n)/(\log 3)<\log(122000)/(\log 3)=10.6605\ldots$, giving that in fact $k\le 15$. Inequality \eqref{eq:low} shows that
$$
\frac{n^{1/2}}{6(\log n)^{1/2}}<15,
$$
so $n<93000$. Next let us observe that if $n$ is not a multiple of $3$, then relation \eqref{eq:mp} leads easily to the conclusion that the number of prime factors $p$ of $C_n$ with $m_p>1$ is in fact $\le (\log n)/(\log 5)=7.15338\ldots$. Hence, the number of such primes is $\le 7$, giving that $k\le 12$, which contradicts a result of Cohen and Hagis \cite{CH} who showed that every number with the Lehmer property must have at least $14$ distinct prime factors.  Hence, $3\mid n$, which shows that $C_n$ is not a multiple of $3$. An argument similar to one used before proves that $n$ is not a multiple of any prime $q>3$. Indeed, for if it were, then relation \eqref{eq:mp} would lead to the conclusion that the number of prime factors $p$ of $C_n$ with $m_p>1$ is $\le 1+\log(n/q)/(\log 3)\le 1+\log(93000/5)/(\log 3)=9.94849\ldots$, so there are at most $9$ such primes. Also, $C_n$ can be divisible with at most $4$ of the $5$ Fermat primes $F_{\gamma}$ with $\gamma\in \{0,1,2,3,4\}$, because $3=F_0$ does not divide $C_n$. Hence, $k\le 9+4=13$, which again contradicts the result from \cite{CH}. Thus, $n=2^{\alpha} 3^{\beta}$ and so all prime factors $p$ of $C_n$ are of the form $2^{\alpha_1}3^{\beta_1}+1$ for some nonnegative integers $\alpha_1$ and $\beta_1$.
Now  write
\begin{equation}
\label{eq:a}
a=\frac{C_n-1}{\phi(C_n)}=\prod_{i=1}^k \left(1+\frac{1}{p_i-1}\right)
\end{equation}
for some integer $a\ge 2$. Since
$$
\prod_{\substack{\alpha_1\ge 0,~\beta_1\ge 0\\ 2^{\alpha_1}3^{\beta_1}+1~{\text{\rm prime}}}}\left(1+\frac{1}{2^{\alpha_1} 3^{\beta_1}}\right)<1.46,
$$
we get that $a<2$, which is a contradiction. This shows that in fact there are no numbers $C_n$ with the claimed property. 

\medskip

We end with some challenges for the reader.

\medskip

{\bf Reearch problem.} {\it Prove that $C_n$ is not a Camichael number for any $n\ge 1$.}

\medskip

If this is too hard, can one at least give a sharp upper bound on the counting function of the set ${\mathcal C}$ of positive integers $n$ such that $C_n$ is a  Carmichael number? We recall that Heppner \cite{Hep} proved that if $x$ is large then the number of positive integers $n\le x$ such that $C_n$ is prime is $O(x/\log x)$, whereas in \cite{LuSh} it was shown that if $a>1$ is a fixed integer then the number of positive integers $n\le x$ such that $C_n$ is base $a$-pseudoprime is $O(x(\log\log x)/\log x)$. Clearly, imposing that $C_n$ is Carmichael (which is a stronger condition) should lead to sharper upper bounds for the counting function of such indices $n$.

\medskip

Finally, here is a problem suggested to us by the referee. Theorem \ref{thm:main} shows that $\phi(C_n)/\gcd(C_n-1,\phi(C_n))$ exceeds $1$ for all $n$. Can one say something more about this ratio? For example, it is possible that a minor modification  of the arguments in the paper would show that this function tends to infinity with $n$, but we have not worked out the details of such a deduction. It would be interesting to find a good (large) lower bound on this quantity which is valid for all $n$ and which tends to infinity with $n$. How about for most $n$? What about lower and upper bounds on the average value of this function when $n$ ranges in the interval $[1,x]$ and $x$ is a large real number? We leave these questions for further research.

\medskip

{\bf Acknowledgements.}  We thank the referee for a careful reading of the paper and for suggesting some of the questions mentioned at the end. F.~L. was supported in part 
by grants PAPIIT 100508 and SEP-CONACyT 79685,


\begin{thebibliography}{9999}

\bibitem{BaGu} W.~D.~Banks, A.~M.~G\"uloglu and C.~W.~Nevans, `On the congruence $n\equiv a\pmod {\varphi(n)}$', {\it INTEGERS} {\bf 8} (2008), \#A59.

\bibitem{BaLu} W.~D.~Banks and F.~Luca, `Composite integers $n$ for which $\varphi(n)\mid n-1$', {\it Acta Math. Sinica\/} {\bf 23} (2007), 1915--1918.

\bibitem{BF} P.~ Berrizbeitia and J.~ G.~ Fernandes, `Observaciones sobre la primalidad de los
n\'umeros de Cullen', Short communication in ``Terceras Jornadas de Teor\'{\i}a de N\'umeros"
(http://campus.usal.es/ tjtn2009/doc/abstracts.pdf).

\bibitem{CiLu} J.~Cilleruelo and F.~Luca, `Repunit Lehmer numbers', {\it Proc. Edinburgh Math. Soc.\/}, to appear.

\bibitem{CH} G.~ L.~ Cohen and P. ~Hagis, `On the number of prime factors of $n$ if
$\phi(n)\mid n-1$', {\it Nieuw Arch. Wisk.\/} {\bf 28} (1980), 177--185.

\bibitem{GO} J.~M.~Grau and A.~M.~Oller-Marc\'en, `An ${\widetilde{O}}(\log^2 N)$  time primality test for generalized Cullen numbers', {\it Preprint\/} 2010, to appear in {\it Math. Comp.\/}

\bibitem{Hep} F. Heppner, `\"Uber Primzahlen der Form $n2^n + 1$ bzw. $p2^p + 1$', {\it Monatsh. 
Math.\/} {\bf 85} (1978), 99--103. 


\bibitem{Ho} C. Hooley, {\it Applications of sieve methods to the theory of numbers\/}, Cambridge University
Press, Cambridge, 1976. 

\bibitem{Lu1} F. Luca, `Fibonacci numbers with the Lehmer property', {\it Bull. Pol. Acad. Sci. Math.\/} {\bf 55} (2007), 7-15.

\bibitem{Lu} F.~Luca, `On the greatest common divisor of two Cullen numbers', {\it Abh. Math. Sem. Univ. Hamburg\/} {\bf 73} (2003), 253--270.

\bibitem{LuPo} F.~Luca and C.~Pomerance, `On composite integers $n$ for which $\phi(n)\mid n-1$", {\it Preprint\/}, 2009.

\bibitem{LuSh} F.~Luca and I.~ E.~ Shparlinski, `Pseudoprime Cullen and Woodall numbers', {\it Colloq.
Math.\/} {\bf 107} (2007), 35--43.

\bibitem{Po1} C.~ Pomerance, `On the distribution of amicable numbers', {\it J. reine angew. Math.\/} {\bf 293/294} (1977), 217--222. 

\bibitem{Po2} C.~Pomerance, `On composite $n$ for which $\varphi(n)\mid n-1$, II', {\it Pacific J. Math.\/} {\bf 69} (1977), 177--186.


\end{thebibliography}
\end{document}